\documentclass[a4paper,oneside,11pt]{article}
\usepackage{amssymb}
\usepackage{amsmath}
\usepackage{array}
\usepackage[utf8]{inputenc}
\usepackage[T1]{fontenc}
\usepackage{enumerate}
\usepackage{graphicx}
\usepackage{theorem}
\usepackage{color}
\usepackage{esint}
\usepackage{accents}

\newtheorem{teorema}{Theorem}

\newtheorem{definition}{Definition}

\newtheorem{lema}{Lemma}
\newtheorem{rema}{Remark}

\newtheorem{proposition}{Proposition}

\newcolumntype{L}{>{\scriptstyle}l}
\newcolumntype{C}{>{\scriptstyle}c}
\newcolumntype{R}{>{\scriptstyle}r}

\def\R{{\mathbb R}}

\def\Z{{\mathbb Z}}

\def\virgp{\raise 2pt\hbox{,}}


\def\XXint#1#2#3{{\setbox0=\hbox{$#1{#2#3}{\int}$ }
\vcenter{\hbox{$#2#3$ }}\kern-.6\wd0}}

\newcommand{\Addresses}{{
  \bigskip
  \footnotesize

 \textsc{IMPA, Rio de Janeiro}. Post-doctoral grant financed by CAPES-Brazil. Partially supported by MICINN grant MTM2013-40945-P (Spain) and ERC grant 277778.\\
  \textit{E-mail address}: \texttt{j.ramosmaravall@gmail.com}

}}

\DeclareMathOperator{\supp}{supp}
\def\<{\langle}
\def\>{\rangle}

\def\({\left(}
\def\){\right)}
\def\supp{\operatorname{supp}}

\makeatletter

\newcommand{\Rmnum}[1]{\expandafter\@slowromancap\romannumeral #1@}
\makeatother

\author{Javier Ramos}
\title{\large \textbf{THE TRILINEAR RESTRICTION ESTIMATE WITH SHARP DEPENDENCE ON THE TRANSVERSALITY} }
\date{}

\begin{document}

\maketitle
\begin{abstract}
We improve the Bennett--Carbery--Tao trilinear restriction estimate for subsets of the paraboloid in three dimensions, giving the sharp factor depending on the transversality. 
\end{abstract}

\section{Introduction}

Let $S\subset \mathbb{R}^3$ be a smooth, compact surface with positive definite second fundamental form, and endowed with its canonical measure $d\sigma$. We write $\mathfrak{R}^*$ for the extension operator,
\begin{align*}
\mathfrak{R}^* f (x)=\int e^{i(\xi,\tau)\cdot x} f(\xi) d\sigma (\xi,\tau),
\end{align*}
acting over functions $f\in \mathcal{S}(\mathbb{R}^2)$. Stein's restriction conjecture asserts that
\begin{align*}
\|\mathfrak{R}^* f\|_{L^{q}(\mathbb{R}^d)}\leq C \|f\|_{L^{\infty}(\R^{d-1})}\quad \mbox{for every}\;\;q> \frac{2d}{d-1}.
\end{align*}
The conjecture with $d=2$ was resolved by Fefferman~\cite{fefferman} (see also \cite{zygmund}). In higher dimensions it remains open, but several partial results have been obtained: Tomas \cite{T}, Stein \cite{Stein2}, Bourgain \cite{bo}, Wolff \cite{wo2}, Moyua--Vargas--Vega \cite{mvv}, \cite{mvv2}, Tao--Vargas--Vega \cite{tvv:bilinear}, Tao--Vargas \cite{tv:cone1}, \cite{tv:cone2}, Tao \cite{taopar}, Bourgain--Guth \cite{boguth} and Guth \cite{G3}. For references regarding the conical and the indefinite case, see \cite{Strichartz}, \cite{Barcelo}, \cite{bo0}, \cite{borg:cone}, \cite{wo}, \cite{sl2} and \cite{V1}. We suggest \cite{T0} for a nice introduction on the subject.

Some of these results were achieved thanks to bilinear techniques. In \cite{taopar}, Tao proved the sharp bilinear estimate, so improvements employing these techniques were exhausted. In \cite{Be}, Bennett introduced the multilinear estimates, and in \cite{bct}, Bennett--Carbery--Tao proved the sharp, up to the endpoint, estimate. In three dimensions it reads as follows:

\Addresses
\clearpage
\begin{teorema}[Bennett-Carbery-Tao]\label{bct1}
Let $S_1,S_2, S_3$ be smooth codimension-one submanifolds of $\mathbb{R}^3$ such that
 \begin{align}\label{H1}
\theta\lesssim |n(\xi_1) \wedge n(\xi_2)  \wedge n(\xi_3) |
\end{align}
for all choices $\xi_n\in S_n$, where $n(\xi_n)$ is the unit normal vector to $S_n$ in $\xi_n$. Then there exist constants $C_{\theta}$ and $\kappa$ such that
\begin{align*}
\|\prod_{n=1}^3 \mathfrak{R}^* f_n \|_{L^{1}(B_R)}\leq  C_{\theta} (\log_2 R)^\kappa \prod_{n=1}^3 \|f_n\|_{L^{2}}
\,\end{align*}
for all $R> 0$.
\end{teorema}
It was unclear for some time how to use this result to deduce linear restriction estimates. That was achieved by Bourgain--Guth in \cite{boguth}. Their argument relies on a dichotomy between good transversality ($\theta\sim 1$) and good $L^4$ orthogonality. Unfortunately, in order to establish the dichotomy there are some inefficiencies which do not allow the full conjecture to be solved.

The main result of this paper is an improvement of the above theorem in dimension $d=3$ with a sharp dependence of the transversality condition \eqref{H1}. Multilinear estimates with this type of transversality dependence were previously considered in related problems, see for example \cite{BHT}, \cite{bb}, \cite{bcct} or \cite{BCW}. For simplicity we carry out the case when \begin{align*}S:=\{(\xi,\tfrac{1}{2}|\xi|^2):\quad \xi\in [0,2]^2\},\end{align*}
but the arguments can be generalized to the case of smooth, compact surfaces with definite second fundamental form. 


\begin{teorema}\label{MT}
Let $S_1,S_2,S_3\subset S$ satisfy
\begin{align}\label{H2}
|n(\xi_1) \wedge n(\xi_2)\wedge n(\xi_3) |\sim\theta
\end{align}
for all choices $\xi_n\in S_n$, where $n(\xi_n)$ is the normal vector to $S_n$ in $\xi_n$. Then there exist constants $C$ and $\kappa$ such that \begin{align}\label{BCTM}
\|\prod_{n=1}^3 \mathfrak{R}^* f_n \|_{L^{1}(B_R)}\leq \theta^{-\frac{1}{2}} C (\log_2 R)^\kappa  \prod_{n=1}^3 \|f_n\|_{L^{2}}
\,\end{align}
for all $R>0$.
\end{teorema}
The factor $\theta^{-\frac{1}{2}}$ is sharp, see Remark \ref{remarkt}. We will apply Theorem \ref{MT} to the linear restriction problem in a forthcoming paper.

The result of Bennett--Carbery--Tao was deduced from a multilinear Kakeya estimate. Later on, Guth \cite{G} (see also Carbery--Valdimarsson \cite{cv} and Guth \cite{G3}) proved the following 

\begin{teorema}[Guth] \label{GTe}
Let $\{T_{n,i}\} $ be a collection of tubes of dimensions $2^\lambda \times 2^\lambda \times 2^{2\lambda}$. Let $v_{n,i}$ be a unit vector parallel to the core of $T_{n,i}$. We assume that the determinant of any of matrices $(v_{1,i_1},v_{2,i_2},v_{3,i_3})$ has norm at least $\theta>0$. Then
\begin{align}\label{MKG}
\int_{\R^3} \prod_{n=1}^3\Big(\sum_{T_{n,i}} \chi_{T_{n,i}}*\mu_{T_{n,i}}\Big)^{\frac{1}{2}}\lesssim \theta^{-\frac{1}{2}} 2^{3\lambda} \prod_{n=1}^3 \Big(\sum_{T_{n,i}} \|\mu_{T_{n,i}}\|\Big)^{\frac{1}{2}}
\end{align}
for all finite measure $\mu_{T_{n,i}}$ and $\lambda\geq 1$.
\end{teorema}

Again, the multilinear Kakeya estimate of Bennett--Carbery--Tao was relevant for the case with $\theta\sim1$. In order to deduce their multilinear restriction estimate, they used an induction on scales argument. Roughly speaking, letting $\mathcal{R}_1(\lambda)$ denote the smallest constant $C$ such that
\begin{align*}
\|\prod_{n=1}^3\mathfrak{R}^* f_n\|_{L^{1}(B_{2^\lambda})}\leq C  \prod_{n=1}^3 \|f_n\|_{L^{2}}
\,\end{align*}
holds, and letting $\mathcal{R}_2(\lambda)$ denote the smallest constant $C$ such that
\begin{align*}\
\int_{\R^3} \prod_{n=1}^3\Big(\sum_{T_{t,i}} \chi_{T_{t,i}}*\mu_{T_{t,i}}\Big)^{\frac{1}{2}}\leq C  2^{3\lambda} \prod_{n=1}^3 \Big(\sum_{T_{t,i}} \|\mu_{T_{t,i}}\|\Big)^{\frac{1}{2}}
\end{align*}
holds, they proved that $\mathcal{R}_1(2\lambda)\lesssim \mathcal{R}_1(\lambda)\mathcal{R}_2(\lambda)$. As by \eqref{MKG}, $\mathcal{R}_2(\lambda)\lesssim 1$, and $\mathcal{R}_1(1)\lesssim 1$, iterating the process they obtained the result. The same argument would not work to obtain Theorem \ref{MT} from \eqref{MKG} as we would be gaining a factor $\theta^{-\frac{1}{2}}$ in each iteration. 

The idea to overcome this problem is linking the hypothesis \eqref{H2} with a refined $L^4$ orthogonality which determines a set in which we iterate the scale. 
More precisely, the strategy is as follows: we decompose $S_1,S_2,S_3$ in subsets for which each triple determines a parallelepiped $P$ in which we have good $L^4$ orthogonality. This orthogonality gives
\begin{align*}
\|\prod_{n=1}^3\mathfrak{R}^* f_n\|_{L^{1}(P)}\lesssim \theta^{-\frac{1}{2}}  \prod_{n=1}^3 \|f_n\|_{L^{2}}.
\,\end{align*}
Now, let $\mathcal{K}(\lambda)$ denote  the smallest constant $C$ such that
\begin{align*}
\|\prod_{n=1}^3 \mathfrak{R}^* f_n\|_{L^{1}(2^\lambda P)}\leq C \theta^{-\frac{1}{2}}    \prod_{n=1}^3 \|f_n\|_{L^{2}}
\,\end{align*}
holds, where $2^\lambda P$ is a dilation of $P$ by $2^\lambda$ . It is enough to prove that $\mathcal{K}(2\lambda)\lesssim \mathcal{K}(\lambda)$. This is accomplished by invoking the Kakeya multilinear estimate \eqref{MKG} together with a discrete version of the multilinear estimate over parallelepiped $2^\lambda P$, which works well because  $|\prod_{n=1}^3\mathfrak{R}^* f_n|$ ``averages'' with no extra factor $\theta^{-\frac{1}{2}}$ in the parallelepiped $2^\lambda P$, as opposed to its average in balls (see Lemma \ref{lemaavO}).

\section{Notation}

Through the paper we will be using the following notation:


$B_r(a)$ is the cube in dimension $3$ of side $r$ and centered in $a$.

 $\tau^j_{k}$ is the square with length side $2^{-j}$ whose left-down vertex is placed in the point $k$. 

 $\tilde{\tau}^j_{k}$ is the lift of $\tau^j_{k}$ to the paraboloid, that is,  $\tilde{\tau}^j_{k}=\{(\xi,\tfrac{1}{2}|\xi|^2):\; \xi\in \tau^j_{k}\}$. 


$t^{j}_{w,m}$ is the strip in the plane of width $2^{-j}$ which passes through $m\in [0,2]^2$ and in the direction $w$, that is
$$t^{j}_{w,m}:=\big\{\xi\in [0,2]^2: \big|\xi-m-w\,((\xi-m)\cdot w)\big|\leq 2^{-j}\big\}.$$
$\tilde{t}^{j}_{w,m}$ is the lift of $t^{j}_{w,m}$ to the paraboloid, that is,  $\tilde{t}^{j}_{w,m}=\{(\xi,\tfrac{1}{2}|\xi|^2):\; \xi\in t^{j}_{w,m}\}$. 

 $\mathbb{S}^1_t$ is a $2^{-t}$ separated set of points on the circumference.

$\varrho^{j}_{w,m}$ is the $2^{-j}$- sector centered in $m$ in the direction $w$, that is


$$\varrho^{j}_{w,m}:=\Big\{\xi\in [0,2]^2: \Big|\frac{\xi-m}{|\xi-m|}-w\Big|\leq 2^{-j}\Big\}.$$



A parallelepiped $\mathcal{P}$ with spanning vectors $(e_1,e_2,e_3)$ and respective length sides ($2^{s_1}$, $2^{s_2}$, $2^{s_3}$), we call a $\Big((s_1,s_2,s_3),(e_1,e_2,e_3)\Big)$ parallelepiped. The dual $\mathcal{P}^*$ is the $\Big((-s_1$, $-s_2$, $-s_3)$, $(e_1$, $e_2$, $e_3)\Big)$ parallelepiped.

We denote the case $s_1=j+2t$, $s_2=j+t$, $s_3=2(j+t)$, $e_1=(w,0)$, $e_2=\big((w,0)\times(m,-1)\big)$, $e_3=(m,-1)$, by $\mathfrak{P}(j,t,w,m)$, the case $s_1=r$, $s_2=r+t$, $s_3=2r$, $e_1=(w,0)$, $e_2=\big((w,0)\times(m,-1)\big)$, $e_3=(m,-1)$, by $\mathfrak{p}(r,t,w,m)$, and the case $s_1=j$, $s_2=j$, $s_3=2j$, $e_1=(1,0,0)$, $e_2=(0,1,0)$, $e_3=(m,-1)$ by $\accentset{\circ}{\tau}_m^j$.

For a $\Big((s_1,s_2,s_3),(e_1,e_2,e_3)\Big)$ parallelepiped $\mathcal{P}$, we write $\mathcal{P}[\lambda]$ for the rescaled $\Big((s_1+\lambda,s_2+\lambda,s_3+\lambda),(e_1,e_2,e_3)\Big)$ parallelepiped.
Also we write $\mathcal{P}(a)$ for parallelepiped centered at $a$.

Let $\phi\in\mathcal{S}$ be radial $\phi\geq 0$, $\phi\geq 1$ in $B_1(0)$, and $\supp \widehat{\phi}\in B_1(0)$, and let $A_{\mathcal{P}}$ be the affine transformation which maps $B_1(0)$ into $\mathcal{P}$. We define $$\phi_{\mathcal{P}}(x)= \phi (A_{\mathcal{P}}^{-1}x).$$


\section{$L^4$ orthogonality}

In this section we prove some $L^4$-type orthogonality estimates.

\begin{figure}[t]
\centering
\includegraphics[width=1.0\textwidth]{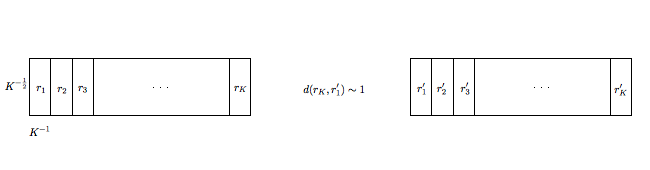}
\caption{The refined orthogonality is represented in the picture: the Minkowski sums of the lifts to the paraboloid $\overline{r}_i+\overline{r}_j'$ are contained in disjoint parallelepipeds. Previously, the known fact was with squares $K^{-\frac{1}{2}}\times K^{-\frac{1}{2}}$.}
\label{fig:awesome_image2}
\end{figure}

\begin{proposition}\label{Port}

. Let $\tau^j_k,\tau^j_{k'}$ be such that $d(\tau^j_k,\tau^j_{k'})=2^{-r}\gtrsim 2^{-j}$, then for every $m\in \tau^j_k,m' \in \tau^j_{k'}$, $w,w'\in \mathbb{S}^1$  with $|w-w'|\lesssim 2^{-t}$, we have
\begin{align*}
&\int\phi_{\mathfrak{P}(j,t,w,m)} \big|\mathfrak{R}^*f_1\chi_{\tau^j_k \cap t^{j+t}_{w,m}} \mathfrak{R}^*f_2\chi_{\tau^j_{k'} \cap t^{j+t}_{w',m}\cap t^{2j-r}_{w'^\perp,m'}}\big|^2\\
&\hspace{10mm}\lesssim \sum_{\substack{\alpha,\;\alpha':\;\alpha\in  2^{-(j+2t)}\mathbb{Z}w\\\alpha'\in  2^{-(2j-r+2t)}\mathbb{Z}w'}} \int\phi_{\mathfrak{P}(j,t,w,m)}  \big|\mathfrak{R}^*f_1\chi_{\tau^j_k \cap t^{j+t}_{w,m}\cap t^{j+2t}_{w^\perp,\alpha}} \mathfrak{R}^*f_2\chi_{\tau^j_{k'}\cap  t^{j+t}_{w',m}\cap t^{2j-r}_{w'^\perp,m'} \cap t^{2j-r+2t}_{w'^\perp,\alpha'} }\big|^2.
\end{align*}
\end{proposition}
For the sake of readability we write the especial case of this proposition when $r=j$ and $w=w'$, see Figure \ref{fig:awesome_image2}.

\begin{proposition}\label{Porto}
Let $\tau^j_k,\tau^j_{k'}$ be such that $2^{-j}= d(\tau^j_k,\tau^j_{k'})$, then for every $m\in \tau^j_k$ and $w\in \mathbb{S}^1$, we have
\begin{align*}
\int\phi_{\mathfrak{P}(j,t,w,m)} &\big|\mathfrak{R}^*f_1\chi_{\tau^j_k \cap t^{j+t}_{w,m}} \mathfrak{R}^*f_2\chi_{\tau^j_{k'} \cap t^{j+t}_{w,m}}\big|^2\\
&\lesssim \sum_{\substack{\alpha,\;\alpha':\;\alpha\in  2^{-(j+2t)}\mathbb{Z}w\\\alpha'\in  2^{-(j+2t)}\mathbb{Z}w}} \int \phi_{\mathfrak{P}(j,t,w,m)}   \big|\mathfrak{R}^*f_1\chi_{\tau^j_k \cap t^{j+t}_{w,m}\cap t^{j+2t}_{w^\perp,\alpha}} \mathfrak{R}^*f_2\chi_{\tau^j_{k'}\cap  t^{j+t}_{w,m}\cap t^{j+2t}_{w^\perp,\alpha'} }\big|^2.
\end{align*}
\end{proposition}

\begin{proposition}\label{Port2}
Let $\tau^j_{k'},\tau^r_{k''}$ be such that $d(\tau^j_{k'},\tau^r_{k''})\sim 2^{-r}$, then for every $w,w'\in \mathbb{S}^1$ with $|w-w'|\lesssim 2^{-t}$, $m \in \tau^j_{k'}$, we have
\begin{align*}
&\int\phi_{\mathfrak{P}(j,t,w,m)} \big|\mathfrak{R}^*f_1\chi_{\tau^j_{k'} \cap t^{j+t}_{w,m}} \mathfrak{R}^*f_2\chi_{\tau^r_{k''} \cap t^{r+t}_{w',m}}\big|^2\\
&\hspace{6mm} \lesssim \sum_{\substack{\alpha,\;w'',\alpha':\;\alpha\in  2^{-(j+2t)}\mathbb{Z}w,\\ w''\in  \mathbb{S}^1_{j+t-r},\\ \alpha'\in   2^{-(2j-r+2t)}\mathbb{Z}w''}} \int \phi_{\mathfrak{P}(j,t,w,m)}  \Big|\mathfrak{R}^*f_1\chi_{\tau^j_{k'} \cap t^{j+t}_{w,m}\cap t^{j+2t}_{w^\perp,\alpha}} \mathfrak{R}^*f_2\chi_{\tau^r_{k''} \cap t^{r+t}_{w',m}\cap \varrho^{j+t-r}_{w'',m}\cap t^{2j-r+2t}_{w''^\perp,\alpha'}}\Big|^2.
\end{align*}
\end{proposition}






\textbf{Proof of Proposition \ref{Port}.}

Using Galilean and rotation invariances we can assume that $m=(0,0)$, $w=(1,0)$ and $\mathfrak{P}(j,t,w,m)$ is a $\Big(j+2t,j+t,2(j+t)\;,\;((1,0,0), (0,1,0), (0,0,1)  )\Big)$ parallelepiped. For $\alpha\in  2^{-(j+2t)}\mathbb{Z}(0,1)$ and $\alpha'\in 2^{-(2j-r+2t)}\mathbb{Z}w'$, let $$\Lambda_\alpha:=\Big\{(\xi,\tfrac{1}{2}|\xi|^2):\xi\in \tau^j_k \cap t^{j+t}_{w,m}\cap t^{j+2t}_{w^\perp,\alpha}\Big\}$$
and 
$$\Lambda_{\alpha'}:=\Big\{(\xi,\tfrac{1}{2}|\xi|^2):\xi\in \tau^j_{k'} \cap t^{j+t}_{w',m}\cap t^{2j-r}_{w'^\perp,m'}\cap t^{2j-r+2t}_{w'^\perp,\alpha'}\Big\}.$$
Observe that by the Galilean and rotation invariances used we have  $\alpha=(\alpha_0,0)$, $\alpha'=\alpha_0' w'$ for some $\alpha_0,\alpha'_0$ with $|\alpha_0|\leq 2^{-j}$, $|\alpha_0'|\sim 2^{-r}$. By Plancherel, the result follows if we prove that for each $\alpha_1,\alpha_1'$
\begin{align}\label{indCo}
\#\Big\{(\alpha_2,\alpha_2'): \quad&\Big(\mathfrak{P}^*(j,t,w,m) +\Lambda_{\alpha_1} + \Lambda_{\alpha_1'} \Big)\bigcap \Big(\mathfrak{P}^*(j,t,w,m) +\Lambda_{\alpha_2} + \Lambda_{\alpha_2'} \Big)\neq \emptyset\Big\}\lesssim 1,
\end{align}
where $+$ is the Minkowski sum.

We begin with the observation that for any $\alpha$, the set $\Lambda_\alpha$ is contained in a mild dilation of some $\Big((-(j+2t),-(j+t),-2(j+t)$), $((1,0,0)$,  $(0,1,0)$,  $(0,0,1)  )\Big)$ parallelepiped, and for any $\alpha'$ the set $\Lambda_{\alpha'}$ is contained in a mild dilation of some  $\Big((-(2j-r+2t)$, $-(j+t)$, $-2(j+t))$, $((w',0)$, $(w'^\perp,0)$, $(0,0,1))\Big)$ parallelepiped. 

Indeed, for any $(\xi,\tfrac{1}{2}|\xi|^2),(\tilde{\xi},\tfrac{1}{2}|\tilde{\xi}|^2)\in \Lambda_{\alpha}$, we have the representation 
\begin{align*}
&\xi=(\alpha_0,0)+(\lambda_1, \lambda_2),\\
&\tilde{\xi}=(\alpha_0,0)+(\lambda_1',\lambda_2'),
\end{align*} 
for some $\lambda_1,\lambda_2,\lambda_1'\lambda_2'$ with $|\lambda_1|,|\lambda_1'|\leq 2^{-(j+2t)}$ and $|\lambda_2|,|\lambda_2'|\leq 2^{-(j+t)}$. 
Therefore, 
\begin{align*}
\big||\xi|^2-|\tilde{\xi}|^2\big|&=\big|\alpha_0 (\lambda_1-\lambda_1')+|\lambda_1|^2+|\lambda_2|^2-|\lambda_1'|^2-|\lambda_2'|^2\big| \lesssim 2^{-2(j+t)}.
\end{align*}

Similarly, for any $(\xi,\tfrac{1}{2}|\xi|^2),(\tilde{\xi},\tfrac{1}{2}|\tilde{\xi}|^2)\in \Lambda_{\alpha'}$ we have the representation 
\begin{align*}
&\xi=\alpha_0'w'+w'\lambda_1+w'^\perp \lambda_2 \\
&\tilde{\xi}=\alpha_0'w'+w'\lambda_1'+w'^\perp \lambda_2' 
\end{align*} 
for some $\lambda_1,\lambda_2,\lambda_1'\lambda_2'$ with $|\lambda_1|,|\lambda_1'|\leq 2^{-(2j-r+2t)}$ and $|\lambda_2|,|\lambda_2'|\leq 2^{-(j+t)}$. Therefore, 
\begin{align*}
\big||\xi|^2-|\tilde{\xi}|^2\big|&=\big|\alpha_0' |\lambda_1-\lambda_1'|+|\lambda_1|^2+|\lambda_2|^2-|\lambda_1'|^2-|\lambda_2'|^2\big| \lesssim 2^{-2(j+t)}.
\end{align*}

We deduce consequently, as by hypothesis $|(1,0)-w'|\leq 2^{-t}$, that the Minkowski sum $\Lambda_\alpha+\Lambda_{\alpha'}$ is also contained in a mild dilation of some $\Big(-(j+2t),-(j+t),-2(j+t)$, $((1,0,0)$,  $(0,1,0)$,  $(0,0,1)  )\Big)$ parallelepiped. 

Now, consider representatives of each subset $(c_{\alpha_n},\tfrac{1}{2}|c_{\alpha_n}|^2)\in \Lambda_{\alpha_n},(c_{\alpha_n'},\tfrac{1}{2}|c_{\alpha_n'}|^2)\in\Lambda_{\alpha_n'}$ for $n=1,2$. In order to prove \eqref{indCo}, after the above observation, it is enough to prove that for any pair $\alpha_1, \alpha_1'$, we just have $O(1)$ indices $\alpha_2, \alpha_2'$ such that
\begin{align}\label{emptys2}\nonumber
&c_{\alpha_1}+c_{\alpha_1'}=c_{\alpha_2}+c_{\alpha_2'}+\Big(O(2^{-(j+2t)}),O(2^{-(j+t)})\Big),\\
&|c_{\alpha_1}|^2+|c_{\alpha_1'}|^2=|c_{\alpha_2}|^2+|c_{\alpha_2'}|^2+O(2^{-2(j+t)}).
\end{align}
We can write $c_{\alpha_2}=c_{\alpha_1}+ \ell (1,0)+O(2^{-(j+t)})(0,1)$, $c_{\alpha_2'}=c_{\alpha_1'}+ \ell'w'+O(2^{-(j+t)})w'^\perp$ for some $\ell,\ell'\in[-2^{-j},2^{-j}]$. It is enough to show $|\ell|\lesssim 2^{-(j+2t)}$ and $|\ell'|\lesssim 2^{-(2j-r+2t)}$. From \eqref{emptys2} we get that $\ell,\ell'$ should obey
\begin{align*}
&\ell=-\ell'+ O( 2^{-(j+2t)}),\\
&\big||\ell|^2+|\ell'|^2+2\ell c_{\alpha_1}(1,0)+2\ell' c_{\alpha_1'}w'+c_{\alpha_1}O(2^{-(j+t)})(0,1)+c_{\alpha_1'}O(2^{-(j+t)})w'^\perp\big|\lesssim 2^{-2(j+t)}.
\end{align*}
As $|\ell|\lesssim 2^{-j}$ and for some $C_1\leq 1, C_2\sim1$ we can write $c_{\alpha_1}=C_12^{-j}(1,0)+O(2^{-(j+t)})(0,1)$ and $c_{\alpha_1'}= (C_2 2^{-r}+C_1 2^{-j})w'+O(2^{-(j+t)})w'^\perp$, we get
\begin{align*}
&\big||\ell'|^2+\ell '(c_{\alpha_1'}w' -c_{\alpha_1}(1,0))\big|\lesssim 2^{-2(j+t)}.
\end{align*}
Then, as $(c_{\alpha_1'}w' -c_{\alpha_1}(1,0))=(c_{\alpha_2'}w' -\ell'-c_{\alpha_1}(1,0))= C 2^{-r}-\ell'$ for some $C\sim 1$, we deduce
\begin{align*}
|\ell' |\lesssim 2^{-2(j+t)}2^{r}.
\end{align*}
From $\ell=-\ell'+ O( 2^{-(j+2t)})$ we also obtain $|\ell|\lesssim 2^{-(j+2t)}$, and the proof is completed.
\begin{flushright}
$\blacksquare$
\end{flushright}

\textbf{Proof of Proposition \ref{Port2}.}
As before, using Galilean and rotation invariances we can assume that $m=(0,0)$ and $w=(1,0)$, and $\mathfrak{P}(j,t,w,m)$ is a $\Big(((j+2t)$, $(j+t)$, $2(j+t))$, $((1,0,0)$, $(0,1,0)$, $(0,0,1)  )\Big)$ parallelepiped.

We first prove the following
\begin{align}\label{FiN}\nonumber
\int\phi_{\mathfrak{P}(j,t,w,m)} &\big|\mathfrak{R}^*f_1\chi_{\tau^j_{k'} \cap t^{j+t}_{w,m}} \mathfrak{R}^*f_2\chi_{\tau^r_{k''} \cap t^{r+t}_{w',m}}\big|^2\\
& \lesssim \sum_{\substack{w'',\alpha':w''\in  \mathbb{S}^1_{j+t-r},\\ \alpha'\in   2^{-(2j-r)}\mathbb{Z}w''}} \int\phi_{\mathfrak{P}(j,t,w,m)}  \Big|\mathfrak{R}^*f_1\chi_{\tau^j_{k'} \cap t^{j+t}_{w,m}} \mathfrak{R}^*f_2\chi_{\tau^r_{k''} \cap t^{r+t}_{w',m}\cap \varrho^{j+t-r}_{w'',m}\cap t^{2j-r}_{w''^\perp,\alpha'}}\Big|^2.
\end{align}
Consider
$$\mathcal{T}:=\big\{(\xi,\tfrac{1}{2} |\xi|^2):\; \xi\in \tau^j_{k'} \cap t^{j+t}_{w,m}\big\}$$
and 
$$\Lambda_{w'',\alpha'}:=\big\{(\xi,\tfrac{1}{2}|\xi|^2):\; \xi\in \tau^r_{k''}\cap t^{r+t}_{w',m}\cap \varrho^{j+t-r}_{w'',m}\cap t^{2j-r}_{w''^\perp,\alpha'}\big\}.$$
According to Plancherel, \eqref{FiN} follows if we prove that for each $w_1'',\alpha_1'$
\begin{align*}
\#\Big\{(w_2'',\alpha_2'): \quad\Big(\mathfrak{P}^*(j,t,w,m) +\mathcal{T} +\Lambda_{w_1'',\alpha_1'} \Big) \bigcap \Big(\mathfrak{P}^*(j,t,w,m) +\mathcal{T}+\Lambda_{w_2'',\alpha_2'} \Big)\neq \emptyset\Big\}\lesssim 1. 
\end{align*}

It is clear that the set $\mathcal{T}$ is contained in mild dilation of a $\Big((-j$, $-(j+t)$, $-2j)$, $((1,0,0)$, $(0,1,0)$, $(0,0,1) )\Big)$ parallelepiped. 

Also, arguing as in the previous proposition, the sets $\Lambda_{w'',\alpha'}$ are contained in a mild dilation of a $\Big((-2j+r,-(j+t),-2j),((w'',0),(w''^{\perp},0),(0,0,1))\Big)$ parallelepiped.

As clearly $|w''-(1,0)|\lesssim 2^{-t}$, we have that the Minkowski sum $\mathcal{T}+\Lambda_{w'',\alpha'}$ is contained in a mild dilation of some $\Big((-j,-(j+t),-2j),((1,0,0), (0,1,0), (0,0,1)\Big)$ parallelepiped.

Now, consider representatives of each subset $(c_{w_n'',\alpha_n'},\tfrac{1}{2}|c_{w_n'',\alpha_n'}|^2)\in \Lambda_{w_n'',\alpha_n'}$ for $n=1,2$. In order to prove \eqref{FiN}, after the above observation, it is enough to prove that for any pair $w_1'', \alpha_1'$, we just have $O(1)$ indices $w_2'', \alpha_2'$ such that
\begin{align}\label{vb1}\nonumber
&c_{w_1'',\alpha_1'}=c_{w_2'',\alpha_2'}+\Big(O(2^{-j}), O(2^{-(j+t)})\Big),\\
&|c_{w_1'',\alpha_1'}|^2=|c_{w_2'',\alpha_2'}|^2+O(2^{-2j}).
\end{align} 
We can write $c_{w_2'',\alpha_2'}=c_{w_1'',\alpha_1'}+ \ell_1w_1''+\ell_2w_1''^\perp$.  It is enough to prove that $|\ell_1|\leq 2^{-(2j-r)}$ and $|\ell_2|\leq 2^{-(j+t)}$. From \eqref{vb1}, we get that 
\begin{align*}
&|\ell_1|\lesssim 2^{-j},|\ell_2|\lesssim 2^{-(j+t)}, \\
&|\ell_1|^2+|\ell_2|^2+2\ell_1c_{w_1'',\alpha_1'}w_1''+2\ell_2c_{w_1'',\alpha_1'}w_1''^\perp\lesssim 2^{-2j}.
\end{align*}
As $c_{w_1'',\alpha_1'}= C 2^{-r}w_1''+O(2^{-(j+t)})w_1''^\perp$ for some $C\sim 1$, we get $|\ell_1|\leq 2^{-(2j-r)}$ and $|\ell_2|\leq 2^{-(j+t)}$, and therefore \eqref{FiN}.

Now, noting that $\varrho^{j+t-r}_{w'',m}\cap \tau^r_{k''}\cap t^{2j-r}_{w''^\perp,\alpha'}$ is contained in a mild dilation of $t^{j+t}_{w'',m}\cap \tau^r_{k''}\cap  t^{2j-r}_{w''^\perp,\alpha'}$, if we apply Proposition \ref{Port}, we infer that for every $w'',\alpha'$,
\begin{align*}
&\int\phi_{\mathfrak{P}(j,t,w,m)}   \Big|\mathfrak{R}^*f_1\chi_{\tau^j_{k'} \cap t^{j+t}_{w,m}} \mathfrak{R}^*f_2\chi_{\tau^r_{k''} \cap t^{r+t}_{w',m}\cap \varrho^{j+t-r}_{w'',m}\cap t^{2j-r}_{w''^\perp,\alpha'}}\Big|^2\\
&\hspace{6mm} \lesssim \sum_{\substack{\alpha,\alpha':\;\alpha\in  2^{-(j+2t)}\mathbb{Z}w,\\ \alpha'\in   2^{-(2j-r+2t)}\mathbb{Z}w''}} \int\phi_{\mathfrak{P}(j,t,w,m)} \Big|\mathfrak{R}^*f_1\chi_{\tau^j_{k'} \cap t^{j+t}_{w,m}\cap t^{j+2t}_{w^\perp,\alpha}} \mathfrak{R}^*f_2\chi_{\tau^r_{k''} \cap t^{r+t}_{w',m}\cap \varrho^{j+t-r}_{w'',m}\cap t^{2j-r+2t}_{w''^\perp,\alpha'}}\Big|^2,
\end{align*}
and therefore the result follows.
\begin{flushright}
$\blacksquare$
\end{flushright}

\section{Trilinear estimate}


\begin{definition}
Let $S_1,S_2,S_3\subset S$ and \eqref{H2} holds. We write $(S_1,S_2,S_3)\sim (r,j,t,w,m,\theta)$ if we can find $(r,j,t,w,m)\in \mathbb{N}\times \mathbb{N}\times \mathbb{N}\times \mathbb{S}^1\times \mathbb{R}^2$, $r\leq j$, such that, perhaps reordering the $S_n$, we have
\begin{align*}
&S_1\subset \{(\xi,\tfrac{1}{2}|\xi|^2):\;\xi\in\tau^{j}_k \cap t^{j+t}_{w,m}\},\\
&S_2\subset \{(\xi,\tfrac{1}{2}|\xi|^2):\;\xi\in\tau^{j}_{k'} \cap t^{j+t}_{w,m}\},\\
&S_3\subset\{(\xi,\tfrac{1}{2}|\xi|^2):\;\xi\in\tau^{r}_{k''} \cap t^{r+t}_{w',m},\},
\end{align*}
for some $k,k',k'',w'$ such that $m\in \tau^{j}_k$, $  d(\tau^{j}_k,\tau^{j}_{k'})\sim 2^{-j}$, $ d(\tau^{j}_{k'},\tau^{r}_{k''})\sim d(\tau^{j}_k,\tau^{r}_{k''})\sim 2^{-r}$, $ |w-w'| \sim 2^{-t}$ and 
\begin{align*}
\theta\sim 2^{-j} 2^{-r} 2^{-t}.
\end{align*}

\end{definition}

\begin{definition}
Let $S_n=	\{(\xi,\tfrac{1}{2}|\xi|^2):\;\xi\in\supp f_n\}$ such that \eqref{H2} holds. We say that $f_1, f_2,f_3$ are $(j,r,t,w,m,\theta)$- triangle type if $(S_1,S_2,S_3)\sim (r,j,t,w,m,\theta)$ for some $(r,j,t,w,m)\in \mathbb{N}\times \mathbb{N}\times \mathbb{N}\times \mathbb{S}^1\times \mathbb{R}^2$.
\end{definition}
The following lemma finds a useful way to write any triple $f_1, f_2,f_3$ whose supports satisfy \eqref{H2} as a sum of triangle type functions.

\begin{lema}\label{3c}
Let $S_n=	\{(\xi,\frac{1}{2}|\xi|^2):\;\xi\in\supp f_n\}$ and \eqref{H2} holds. Then, there exists a collection $\{S_{1,i},S_{2,i},S_{3,i},r_i,j_i,t_i,w_i,m_i\}_{i=1}$ such that 

i) It is a partition
\begin{align*}
S_1\times S_2 \times S_3=\bigcup_{i\leq C(\theta)} S_{1,i}\times S_{2,i} \times S_{3,i}.
\end{align*}

ii) $(S_{1,i},S_{2,i},S_{3,i})\sim (r_i,j_i,t_i,w_i,m_i,\theta)$.


iii) There exist subsets of indices $\{I_j\}_{j\leq C\log_2 \theta^{-1}}$ such that $$\bigcup_{i\leq C(\theta)} S_{1,i}\times S_{2,i} \times S_{3,i}= \bigcup_{j\leq C\log_2 \theta^{-1}} \bigcup_{i\in I_j} S_{1,i}\times S_{2,i} \times S_{3,i}$$
and for every $j\leq C\log_2 \theta^{-1}$, $i\in I_j$,
\begin{align}\label{losns}\nonumber
&\#\{i'\in I_j:\;\; S_{n_1,i}\cap S_{n_1,i'}\neq \emptyset\}\lesssim \log_2 \theta^{-1},\\
&\#\{i'\in I_j:\;\; S_{n_2,i}\cap S_{n_2,i'}\neq \emptyset\} \lesssim \log_2 \theta^{-1}\quad \mbox{for some}\;\;n_1,n_2\in \{1,2,3\}\;\;n_1\neq n_2.
\end{align}

\end{lema}
\begin{figure}[t]\label{figS}
\centering
\includegraphics[width=1.0\textwidth]{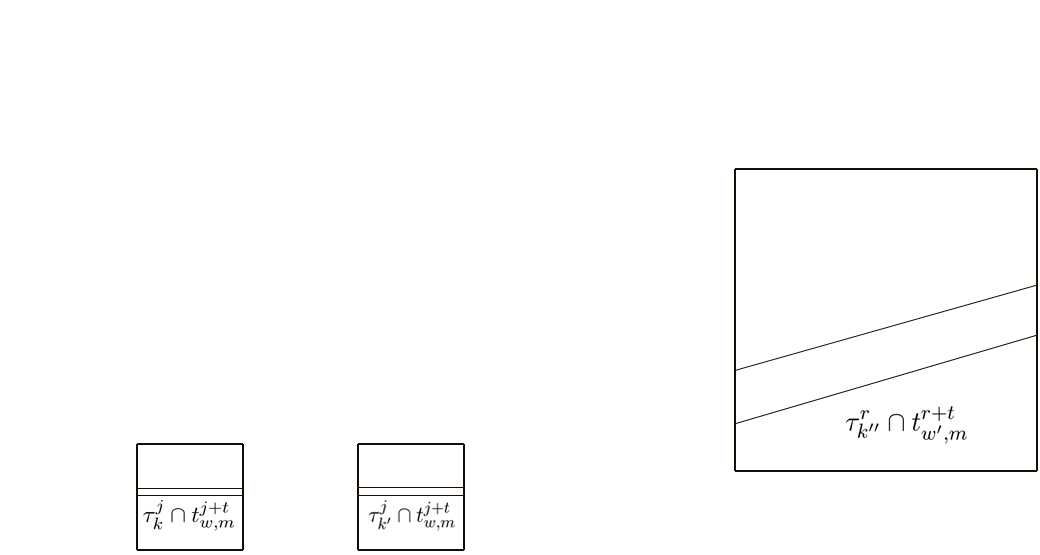}
\caption{Example of $S_1, S_2, S_3$ with $(S_1,S_2,S_3)\sim (r,j,t,w,m,\theta)$.}
\end{figure}
\textbf{Proof of Lemma \ref{3c}}.

It is easy to check that the hypothesis \eqref{H2} means that for every $\xi_1\in \supp f_1$, $\xi_2\in \supp f_2$ and $\xi_3\in \supp f_3$, the area $T$ of the triangle with vertices $\xi_1,\xi_2,\xi_3$ is $T\sim\theta$.

We use a Whitney type decomposition (see \cite{tvv:bilinear}) adapted to the trilinear setting.
We write $\tau_{k}^j\sim \tau_{k'}^j$ if $\tau_{k}^j$ is not adjacent to $\tau_{k'}^j$ but have adjacent parents. For almost every $x, y\in [0,2]^2$, there exists a unique pair $\tau_{k}^j, \tau_{k'}^j$ with $\tau_{k}^j\sim\tau_{k'}^j$ containing $x$ and $y$ respectively.

For every $j,k=(k_1,k_2),k'=(k_1',k_2')$ with $\tau_{k}^j\sim \tau_{k'}^j$ and $\mathbf{k}_j=(\min(k_1,k_1'),\min(k_2,k_2'))-(2^{-j+2},2^{-j+2})$, we write $\tau^j_{k''}\sim[\tau_{k}^j, \tau_{k'}^j]$ if $\tau^{j}_{k''}$ is contained in the square $\tau^{j-4}_{\mathbf{k}_j}$ but is not adjacent to $\tau_{k}^j$ or $\tau_{k'}^j$, and for $r<j$ recursively we write $\tau^r_{k''}\sim[\tau_{k}^j, \tau_{k'}^j]$ if $\tau_{k''}^r$ is adjacent to some $\tau_{k'''}^{r+1}$ with $\tau_{k'''}^{r+1}\sim[\tau_{k}^j, \tau_{k'}^j]$, it does not contain any $\tau_{k''''}^{r+1}$ with $\tau_{k''''}^{r+1}\sim[\tau_{k}^j, \tau_{k'}^j]$ and $k''\in 2^{-r}\Z^2+\mathbf{k}_{r+1}$, and we set $\mathbf{k}_{r}=\mathbf{k}_{r+1}-(2^{-r},2^{-r})$.






We define
\begin{align*}
&A_1:=\bigcup_{j\geq 0}\bigcup_{k}\;\;\bigcup_{k':\tau^j_k\sim\tau^j_{k'}} \;\bigcup_{0\leq r\leq j}\;\;\bigcup_{k'': \tau^r_{k''}\sim[\tau_{k}^j, \tau_{k'}^j]}\tau^j_k\times\tau^j_{k'}\times \tau^r_{k''},\\
&A_2:=\bigcup_{j\geq 0}\bigcup_{k}\;\;\bigcup_{k':\tau^j_k\sim\tau^j_{k'}} \;\bigcup_{0\leq r\leq j}\;\;\bigcup_{\substack{k'': \tau^r_{k''}\sim[\tau_{k}^j, \tau_{k'}^j]\\ \mbox{\footnotesize{if}}\;r=j,\;\mbox{\footnotesize{then}}\; \tau^j_{k'}\not\sim[\tau_{k}^j, \tau_{k''}^r]}}\tau^j_k\times\tau^r_{k''}\times \tau^j_{k'},\\
&A_3:=\bigcup_{j\geq 0}\bigcup_{k}\;\;\bigcup_{k':\tau^j_k\sim\tau^j_{k'}} \;\bigcup_{0\leq r\leq j}\;\;\bigcup_{\substack{k'': \tau^r_{k''}\sim[\tau_{k}^j, \tau_{k'}^j] \\\mbox{\footnotesize{if}}\;r=j,\;\mbox{\footnotesize{then}}\,\tau^j_{k'}\not\sim[\tau_{k}^j, \tau_{k''}^r]\;\mbox{\footnotesize{and}}\;\tau^j_{k}\not\sim[\tau_{k'}^j, \tau_{k''}^r]}}\tau^r_{k''}\times\tau^j_k\times \tau^j_{k'}.
\end{align*}






We have the following partition
$$[0,2]^6=A_1\bigcup A_2\bigcup A_3\quad \mbox{almost everywhere}.$$

Indeed, by construction we have $$[0,2]^6\setminus A_1= \bigcup_{j\geq 0}\bigcup_{k}\;\;\bigcup_{k':\tau^j_k\sim\tau^j_{k'}}\;\;\;\bigcup_{k'':\; \tau^j_{k''}\;\mbox{adjacent to }\;\tau^j_k\;\mbox{or}\;\tau^j_{k'}}\tau^j_{k}\times\tau^j_{k'}\times\tau^j_{k''},$$ 
and for each term $\tau^j_{k}\times\tau^j_{k'}\times\tau^j_{k''}$ in the previous partition, it is easy to check that we have
\begin{align*}
&\tau^j_{k}\times\tau^j_{k'}\times\tau^j_{k''}\\
=&\bigcup_{j'> j}\bigcup_{k'''}\;\;\bigcup_{k'''':\tau^{j'}_{k'''}\sim\tau^{j'}_{k''''}}\bigcup_{r\leq j'}\bigcup_{k''''': \tau_{k'''''}^r\sim [\tau^{j'}_{k'''},\tau^{j'}_{k''''}]}(\tau^j_{k}\cap \tau^r_{k'''''}) \times(\tau^{j}_{k'}\cap \tau^{j'}_{k'''})\times(\tau^{j}_{k''}\cap \tau^{j'}_{k''''})\\
&\bigcup_{j'> j}\bigcup_{k'''}\;\;\bigcup_{k'''':\tau^{j'}_{k'''}\sim\tau^{j'}_{k''''}}\bigcup_{r\leq j'}\bigcup_{k''''': \tau_{k'''''}^r\sim [\tau^{j'}_{k'''},\tau^{j'}_{k''''}]}(\tau^j_{k}\cap \tau^{j'}_{k'''} ) \times(\tau^j_{k'}\cap \tau^r_{k'''''})\times(\tau^j_{k''}\cap \tau^{j'}_{k''''}).
\end{align*}
Thus, $\tau^j_{k}\times\tau^j_{k'}\times\tau^j_{k''}\subset A_2\bigcup A_3$, and we have proven that for almost every $x,y,z\in [0,2]^2$ we have $x\times y\times z \in A_n$ for some $n\in\{1,2,3\}$. It just remains to check that $|A_{n}\cap A_{n'}|=0$ for $n\neq n'$. This follows by observing that for any $j$, $j'$, $r\leq j-1$, $r'$, $k$, $k'$, $k''$, $k'''$, $k''''$, $k'''''$ with $\tau_k^j\sim\tau_{k'}^j$, $\tau_{k''}^r\sim [\tau_k^j,\tau_{k'}^j]$, $\tau_{k'''}^{j'}\sim \tau_{k''''}^{j'}$ and $\tau_{k'''''}^{r'}\sim [\tau_{k'''}^{j'}, \tau_{k''''}^{j'}]$, if $j'\geq j$ then as the parents of $\tau_{k'''}^{j'}$ and $\tau_{k''''}^{j'}$ are adjacent and $d(\tau_k^j,\tau_{k''}^r), d(\tau_{k'}^j,\tau_{k''}^r)\geq 2^{-j+2}$, we get $\tau_{k''''}^{j'}\cap \tau_{k''}^r=\emptyset$ or $\tau_{k'''}^{j'}\cap \tau_{k}^j=\emptyset$ and $\tau_{k'''}^{j'}\cap \tau_{k'}^j=\emptyset$, and if $j'< j$ then as the parents of $\tau_{k}^{j}$ and $\tau_{k'}^{j}$ are adjacent we get $\tau_{k'''''}^{r'}\cap \tau_{k}^j=\emptyset$ or $\tau_{k''''}^{j'}\cap \tau_{k'}^j=\emptyset$ and $\tau_{k'''}^{j'}\cap \tau_{k'}^j=\emptyset$.



Now, for each term $\tau^j_k\times\tau^j_{k'}\times \tau^r_{k''}$ in $A_1$, let $t$ be such that $2^{-t}\sim \theta 2^{j} 2^{r}$, we decompose $\tau^j_k\times\tau^j_{k'}$,
\begin{align*}
\tau^j_k\times\tau^j_{k'}=\bigcup_{w\in \mathbb{S}^1_{t+2},m\in 2^{-(t+j)}\mathbb{Z}w^\perp} \tau^j_k\cap t^{j+t}_{w,m} \times\tau^j_{k'}\cap t^{j+t}_{w,m}.
\end{align*}
As for each $w\in \mathbb{S}^1_{t+2},m\in 2^{-(t+j)}\mathbb{Z}w^\perp$,
\begin{align*}
\#\big\{(w',m'): \; w'\in\mathbb{S}^1_{t+2},\, m'\in 2^{-(t+j)}\mathbb{Z}w'^\perp : &\;\;\tau^j_k\cap t^{j+t}_{w,m}\cap t^{j+t}_{w',m'}\neq \emptyset,\\
&\;\; \tau^j_{k'}\cap t^{j+t}_{w,m}\cap t^{j+t}_{w',m'}\neq \emptyset \big\}\lesssim 1,
\end{align*}
we can break the sets $\tau^j_k\cap t^{j+t}_{w,m} \times\tau^j_{k'}\cap t^{j+t}_{w,m}$ in $O(1)$ subsets and get a partition,
\begin{align*}
\tau^j_k\times\tau^j_{k'}=\bigcup_{w\in \mathbb{S}^1_{t+2},m\in 2^{-(t+j)}\mathbb{Z}w^\perp}\bigcup_{s\leq C} \tau^j_k\cap t^{j+t}_{w,m} \cap E_{j,k,k',w,m,s} \times\tau^j_{k'}\cap t^{j+t}_{w,m}\cap E_{j,k,k',w,m,s}
\end{align*}
for some sets $E_{j,k,k',w,m,s}$. Fix $(w,m)$, by \eqref{H2}, we have 
\begin{align}\label{angle123}
|\angle (\xi_3-\xi_1,\xi_2-\xi_1)|\sim 2^{-t}
\end{align}
for every $\xi_1,\xi_2,\xi_3$ such that $(\xi_1,\frac{1}{2}|\xi_1|^2)\in \tilde{\tau}^j_k\cap \tilde{t}^{j+t}_{w,m} \cap \tilde{E}_{j,k,k',w,m,s}\cap S_1$, $(\xi_2,\frac{1}{2}|\xi_2|^2)\in \tilde{\tau}^j_{k'}\cap \tilde{t}^{j+t}_{w,m}\cap \tilde{E}_{j,k,k',w,m,s}\cap S_2$ and $(\xi_3,\frac{1}{2}|\xi_3|^2)\in \tilde{\tau}^r_{k''}\cap S_3$, where $\tilde{E}_{j,k,k',w,m,s}$ is the lift of $E_{j,k,k',w,m,s}$ to the paraboloid. Therefore, we get 
\begin{align*}
\tilde{\tau}^r_{k''}\cap S_3\subset \tilde{\tau}^r_{k''}\cap \bigcup_{\substack{w'\in  \mathbb{S}^1_{t+2}, \; m' \in 2^{-(t+r)}\mathbb{Z}w'^\perp,\\|w'-w|\sim 2^{-t},\\t^{r+t}_{w',m'}\cap \tau^j_k\cap t^{j+t}_{w,m}\neq \emptyset}}\tilde{t}^{r+t}_{w',m'}\cap E_{j,k,k',w,m,w',m'},
\end{align*}
for some sets $E_{j,k,k',w,m,w',m'}$ such that $\{\tilde{t}^{r+t}_{w',m'}\cap E_{j,k,k',w,m,w',m'}\}_{w',m'}$ are disjoint sets. Clearly, there are just $O(1)$ indices $w',m'$ like that. We associate each $j,k,k',r,k'',w,m,s,w',m'$ with an index $i$, and define
 \begin{align*}
& S_{1,i}= \tilde{\tau}^j_k\cap \tilde{t}^{j+t}_{w,m} \cap \tilde{E}_{j,k,k',w,m,s}\cap S_1,\\
& S_{2,i}=  \tilde{\tau}^j_{k'}\cap \tilde{t}^{j+t}_{w,m}\cap \tilde{E}_{j,k,k',w,m,s}\cap S_2,\\
& S_{3,i}= \tilde{\tau}^r_{k''}\cap\tilde{t}^{r+t}_{w',m'}\cap \tilde{E}_{j,k,k',w,m,w',m'}\cap S_3.
 \end{align*}
 We argue equivalently with $A_2$ and $A_3$. The sets obtained $\{S_{1,i}, S_{2,i}, S_{3,i}\}_i$ satisfy i) and ii). In order to see iii), we first notice that we consider the cases $2^j\lesssim \theta^{-1}$. Also, for each $j,k$, there are $O(1)$ indices $k'$ such that $\tau^j_k\sim \tau^{j}_{k'}$ and $O(j)$ indices $r,k''$ such that $\tau^r_{k''}\sim [\tau^j_k, \tau^{j}_{k'}]$. Finally, for each $j,k,k',r,k''$ there are $O(1)$ indices $w,m$. Indeed, fix $\xi_1,\xi_2,\xi_3$ such that $(\xi_1,\frac{1}{2}|\xi_1|^2)\in\tilde{\tau}^j_{k}\cap S_1$, $(\xi_2,\frac{1}{2}|\xi_2|^2)\in\tilde{\tau}^j_{k'}\cap S_2$ and $(\xi_3,\frac{1}{2}|\xi_3|^2)\in \tilde{\tau}^r_{k''}\cap S_3$, then $\xi_3,\xi_1$ impose that $S_2\subset  \tilde{\tau}^j_{k'}\cap\displaystyle \bigcup_{i\leq C} \tilde{t}^{j+t}_{w_i,m_i}$, and $\xi_3,\xi_2$ impose that $S_1\subset  \tilde{\tau}^j_{k}\cap\displaystyle \bigcup_{i\leq C} \tilde{t}^{j+t}_{w_i,m_i}$. For each $A_{n}$ the subsets of indices $\{I_j\}$ will be those with fixed $j$, and the $n_1,n_2$ in $\eqref{losns}$ will be those from the Whitney decomposition.
\begin{flushright}
$\blacksquare$
\end{flushright}

 \begin{rema}
Using the ideas of the proof of Lemma \ref{3c} we can show that hypothesis \eqref{H2} in Theorem \ref{MT} could be substitute by $$|n(\xi_1) \wedge n(\xi_2)\wedge n(\xi_3) |\gtrsim\theta.$$
 \end{rema}

In our notation, we immediately get the following from Theorem \ref{GTe}.

\begin{teorema}\label{GuthT}
Let $f_1,f_2,f_3$ of $(j,r,t,w,m,\theta)$-triangle type, then
\begin{align*}
\int_{\R^3} \prod_{n=1}^3\Big(\sum_{k\in 2^{-\lambda} \mathbb{Z}^2\cap\supp f_n} \mu_{\accentset{\circ}{\tau}^{\lambda}_{k}}*\chi_{\accentset{\circ}{\tau}^{\lambda}_{k}}\Big)^{\frac{1}{2}}\lesssim 2^{\frac{j+r+t}{2}} 2^{3\lambda} \prod_{n=1}^3 \Big(\sum_{k\in 2^{-\lambda} \mathbb{Z}^2\cap\supp f_n}\|\mu_{\accentset{\circ}{\tau}^{\lambda}_{k}}\|\Big)^{\frac{1}{2}}
\end{align*}
for all finite measure $\mu_{\accentset{\circ}{\tau}^{\lambda}_{k}}$.

\end{teorema}

Theorem \ref{MT} will follow from the following theorem.
\begin{teorema}\label{TriR}
If $S_1,S_2,S_3$ are $(j,r,t,w,m,\theta)$-triangle type then there exist constants $C$ and $\kappa$ such that
\begin{align}\label{multi}
\int_{B_R} |\mathfrak{R}^*f_1\mathfrak{R}^*f_2 \mathfrak{R}^*f_3|\lesssim C (\log_2 R)^\kappa 2^{\frac{j+r+t}{2}}\|f_1\|_{L^2}\|f_2\|_{L^2}\|f_3\|_{L^2}
\end{align}
for all $R>0$.

\end{teorema}

\textbf{Proof of Theorem \ref{MT}.}

First, we notice that if $\theta\leq R^{-10}$ then the result is trivial using $|\mathfrak{R}^*F|(x) \leq \|F\|_1\leq |\supp F|^{\frac{1}{2}}\|F\|_2$.

By Lemma \ref{3c} i), ii) and the triangle inequality
\begin{align*}
\int_{B_R} |\mathfrak{R}^*f_1\mathfrak{R}^*f_2 \mathfrak{R}^*f_3|\leq \sum_{i}\int_{B_R} |\mathfrak{R}^*f_{1,i}\mathfrak{R}^*f_{2,i} \mathfrak{R}^*f_{3,i}|,
\end{align*}
where for each $i$, the functions $f_{1,i},f_{2,i},f_{3,i}$ are of $(j_i,r_i,t_i,w_i,m_i,\theta)$-triangle type for some $(j_i,r_i,t_i,w_i,m_i)$. By Theorem \ref{TriR},
\begin{align*}
\int_{B_R} |\mathfrak{R}^*f_1\mathfrak{R}^*f_2 \mathfrak{R}^*f_3|\lesssim \theta^{-\frac{1}{2}}(\log_2 R)^\kappa \sum_{i} \|f_{1,i}\|_{L^2}\|f_{2,i}\|_{L^2}\|f_{3,i}\|_{L^2}.\end{align*}
By Lemma \ref{3c} iii) and Cauchy-Schwarz inequality we deduce the result.
\begin{flushright}
$\blacksquare$
\end{flushright}


We use the following definition.
\begin{definition}
We denote by $\mathcal{K}(\lambda)$ the smallest constant $C$ such that
\begin{align*}
\int_{\mathfrak{P}(j,t,w,m)[\lambda]}|\mathfrak{R}^*f_1\mathfrak{R}^*f_2 \mathfrak{R}^*f_3|\leq C 2^{\frac{j+r+t}{2}}\|f_1\|_{L^2}\|f_2\|_{L^2}\|f_3\|_{L^2}
\end{align*}
for every $f_1,f_2,f_3$  of $(j,r,t,w,m,\theta)$-triangle type. 

\end{definition}

The induction step explained in the introduction is the next proposition.

\begin{proposition}\label{P1}


$\mathcal{K}(2\lambda)\lesssim \mathcal{K}(\lambda)$.
\end{proposition}

We also need the initial condition to start the induction argument.
\begin{proposition}\label{P2}

$\mathcal{K}(1)\lesssim 1$ 

\end{proposition}

\textbf{Proof of Theorem \ref{TriR}.}

It is enough to prove $\mathcal{K}(\log_2 R)\lesssim (\log_2 R)^\kappa$. By Propositions \ref{P1} applied $O(\log_2\log_2 R)$ times and by Proposition \ref{P2} the result follows.

\begin{flushright}
$\blacksquare$
\end{flushright}

For the proof of Proposition \ref{P2} and in a forthcoming paper we will need the following\footnote{In the present paper we just need $L^1$ average and the `constant' property will not be used.}.
\begin{proposition}\label{TriOrt}
Let $f_1, f_2,f_3$ be $(j,r,t,w,m,\theta)$-triangle type. Then, there exists a function $\psi$, such that 
\begin{align*}
& - L^{\frac{4}{3}}\;\mbox{average on }\;\mathfrak{P}(j,t,w,m):\quad\fint_{\mathfrak{P}(j,t,w,m)}|\psi|^{\frac{4}{3}}\lesssim 1,\\
& - \mbox{`Constant' on }\;\mathfrak{p}(j,t,w,m):\quad \psi=\psi*|\mathfrak{p}(r,t,w,m)|^{-1}\varphi_{\mathfrak{p}(r,t,w,m)}\\
&\hspace{55mm}\mbox{for some Schwartz function }\; \varphi_{\mathfrak{p}(r,t,w,m)} \;\mbox{adapted to}\; \mathfrak{p}(r,t,w,m),
\end{align*}
and 
\begin{align*}
&\big|\mathfrak{R}^*f_1 \mathfrak{R}^*f_2 \mathfrak{R}^*f_3\big|(x)=\psi (x) \sup_{y_1,y_2,y_3}\prod_{n=1}^3\phi_{\mathfrak{P}(j,t,w,m)}(y_n)\\
&\Big(\sum_{\alpha\in2^{-(j+2t)}\mathbb{Z}w}\big|\mathfrak{R}^* f_{1}\chi_{\tau^j_k \cap t^{j+t}_{w,m}\cap t^{j+2t}_{w^\perp,\alpha}}(y_1)\big|^2 \Big)^{\frac{1}{2}}\Big(\sum_{\alpha' \in2^{-(j+2t)}\mathbb{Z}w}\big| \mathfrak{R}^* f_{2}\chi_{\tau^j_{k'} \cap t^{j+t}_{w,m}\cap t^{j+2t}_{w^\perp,\alpha'}}(y_2)\big|^2\Big)^{\frac{1}{2}} \\
& \Big(\sum_{w''\in \mathbb{S}^1_{j+t-r},\; \alpha'' 2^{-(2j-r+2t)}\mathbb{Z}w''}\big|\mathfrak{R}^* f_{3}\chi_{\tau^r_{k''} \cap t^{r+t}_{w',m}\cap \varrho^{j+t-r}_{w'',m}\cap t^{2j-r+2t}_{w''^\perp,\alpha''}}(y_3)\big|^2 \Big)^{\frac{1}{2}}.
\end{align*}
\end{proposition}


\textbf{Proof of Proposition \ref{TriOrt}.}

By definition, there exists $k,k',k'', w'$ such that
\begin{align*}
\fint_{\mathfrak{P}(j,t,w,m)} \big|\mathfrak{R}^*f_1 \mathfrak{R}^*f_2 \mathfrak{R}^*f_3\big|^{\frac{4}{3}}=\fint_{\mathfrak{P}(j,t,w,m)} \big|\mathfrak{R}^*f_1\chi_{\tau^j_k \cap t^{j+t}_{w,m}}  \mathfrak{R}^*f_2\chi_{\tau^j_{k'} \cap t^{j+t}_{w,m}} \mathfrak{R}^*f_3\chi_{\tau^r_{k''} \cap t^{r+t}_{w',m}}\big|^{\frac{4}{3}}.
\end{align*}

By H\"older's inequality,
\begin{align*}
&\fint_{\mathfrak{P}(j,t,w,m)} \big|\mathfrak{R}^*f_1\chi_{\tau^j_k \cap t^{j+t}_{w,m}}  \mathfrak{R}^*f_2\chi_{\tau^j_{k'} \cap t^{j+t}_{w,m}} \mathfrak{R}^*f_3\chi_{\tau^r_{k''} \cap t^{r+t}_{w',m}}\big|^{\frac{4}{3}}\\
\leq &\Big(\fint_{\mathfrak{P}(j,t,w,m)} \big|\mathfrak{R}^*f_1\chi_{\tau^j_k \cap t^{j+t}_{w,m}}  \mathfrak{R}^*f_2\chi_{\tau^j_{k'} \cap t^{j+t}_{w,m}} \big|^{2}\Big)^{\frac{1}{3}}\Big(\fint_{\mathfrak{P}(j,t,w,m)} \big|\mathfrak{R}^*f_1\chi_{\tau^j_k \cap t^{j+t}_{w,m}}  \mathfrak{R}^*f_3\chi_{\tau^r_{k''} \cap t^{r+t}_{w',m}}\big|^{2}\Big)^{\frac{1}{3}}\\
&\hspace{65mm} \Big(\fint_{\mathfrak{P}(j,t,w,m)} \big|\mathfrak{R}^*f_2\chi_{\tau^j_{k'} \cap t^{j+t}_{w,m}} \mathfrak{R}^*f_3\chi_{\tau^r_{k''} \cap t^{r+t}_{w',m}}\big|^{2}\Big)^{\frac{1}{3}}\\
=&\; I_1 \cdot I_2 \cdot I_3.
\end{align*}
For $I_1$ we use Proposition \ref{Porto} to get
\begin{align*}
I_1 \lesssim \Big(\sum_{\alpha,\;\alpha':\;\alpha,\alpha'\in  2^{-(j+2t)}\mathbb{Z}w} \fint\phi_{\mathfrak{P}(j,t,w,m)}  \big|\mathfrak{R}^*f_1\chi_{\tau^j_k \cap t^{j+t}_{w,m}\cap t^{j+2t}_{w^\perp,\alpha}} \mathfrak{R}^*f_2\chi_{\tau^j_{k'} \cap t^{j+t}_{w,m}\cap t^{j+2t}_{w^\perp,\alpha'}}\big|^2\Big)^{\frac{1}{3}}.
\end{align*}
For $I_2$ and $I_3$ we use Proposition \ref{Port2} to get
\begin{align*}
I_2 & \lesssim  \Big( \sum_{\substack{\alpha,\;w'',\alpha'':\;\alpha\in  2^{-(j+2t)}\mathbb{Z}w,\\ w''\in  \mathbb{S}^1_{j+t-r},\\ \alpha''\in   2^{-(2j-r+2t)}\mathbb{Z}w''}} \fint\phi_{\mathfrak{P}(j,t,w,m)}  \Big|\mathfrak{R}^*f_1\chi_{\tau^j_{k} \cap t^{j+t}_{w,m}\cap t^{j+2t}_{w^\perp,\alpha}} \mathfrak{R}^*f_3\chi_{\tau^r_{k''} \cap t^{r+t}_{w',m}\cap \varrho^{j+t-r}_{w'',m}\cap t^{2j-r+2t}_{w''^\perp,\alpha''}}\Big|^2\Big)^{\frac{1}{3}}
\end{align*}
and
\begin{align*}
I_3 & \lesssim  \Big(\sum_{\substack{\alpha',\;w'',\alpha'':\;\alpha'\in  2^{-(j+2t)}\mathbb{Z}w,\\ w''\in  \mathbb{S}^1_{j+t-r},\\ \alpha''\in   2^{-(2j-r+2t)}\mathbb{Z}w''}} \fint\phi_{\mathfrak{P}(j,t,w,m)}  \Big|\mathfrak{R}^*f_2\chi_{\tau^j_{k'} \cap t^{j+t}_{w,m}\cap t^{j+2t}_{w^\perp,\alpha'}}  \mathfrak{R}^*f_3\chi_{\tau^r_{k''} \cap t^{r+t}_{w',m}\cap \varrho^{j+t-r}_{w'',m}\cap t^{2j-r+2t}_{w''^\perp,\alpha''}}\Big|^2\Big)^{\frac{1}{3}}.
\end{align*}
Taking suprema we conclude the result observing that we can find a Schwartz function $\varphi_{\mathfrak{p}(r,t,w,m)}$ adapted to $\mathfrak{p}(r,t,w,m)$ with $$\mathfrak{R}^*f_1 \mathfrak{R}^*f_2 \mathfrak{R}^*f_3(x)=\mathfrak{R}^*f_1 \mathfrak{R}^*f_2 \mathfrak{R}^*f_3*|\mathfrak{p}(r,t,w,m)|^{-1}\varphi_{\mathfrak{p}(r,t,w,m)}(x).$$
\begin{flushright}
$\blacksquare$
\end{flushright}

\textbf{Proof of Proposition \ref{P2}.}

By Proposition \ref{TriOrt}, 
\begin{align*}
\int_{\mathfrak{P}(j,t,w,m)}\big|\mathfrak{R}^*f_1 \mathfrak{R}^*f_2 \mathfrak{R}^*f_3\big|\lesssim \big|\mathfrak{P}(j,t,w,m)\big|\sup_{y_1,y_2,y_3}\big|\Big(\sum_{\alpha} \big|\mathfrak{R}^* f_{1}\chi_{\tau^j_k \cap t^{j+t}_{w,m}\cap t^{j+2t}_{w^\perp,\alpha}}(y_1)\big|^2 \Big)^{\frac{1}{2}}&\\		
\Big(\sum_{\alpha'}\big| \mathfrak{R}^* f_{2}\chi_{\tau^j_{k'} \cap t^{j+t}_{w,m}\cap t^{j+2t}_{w^\perp,\alpha'}}(y_2)\big|^2\Big)^{\frac{1}{2}}&\\ \Big(\sum_{w'',\alpha''} \big|\mathfrak{R}^* f_{3}\chi_{\tau^r_{k''} \cap t^{r+t}_{w',m}\cap \varrho^{j+t-r}_{w'',m}\cap t^{2j-r+2t}_{w''^\perp,\alpha''}}(y_3)\big|^2 \Big)^{\frac{1}{2}}&.
\end{align*}
Using $|\mathfrak{R}^*F|(x) \leq \|F\|_1$, we have
\begin{align*}
\int_{\mathfrak{P}(j,t,w,m)}\big|\mathfrak{R}^*f_1 \mathfrak{R}^*f_2 \mathfrak{R}^*f_3\big|\lesssim \big|\mathfrak{P}(j,t,w,m)\big|\Big(\sum_{\alpha} \| f_{1}\chi_{\tau^j_k \cap t^{j+t}_{w,m}\cap t^{j+2t}_{w^\perp,\alpha}}\|_1^2 \Big)^{\frac{1}{2}}&\\		
\Big(\sum_{\alpha'}\| f_{2}\chi_{\tau^j_{k'} \cap t^{j+t}_{w,m}\cap t^{j+2t}_{w^\perp,\alpha'}}\|_1^2\Big)^{\frac{1}{2}}&\\
  \Big(\sum_{w'',\alpha''} \| f_{3}\chi_{\tau^r_{k''} \cap t^{r+t}_{w',m}\cap \varrho^{j+t-r}_{w'',m}\cap t^{2j-r+2t}_{w''^\perp,\alpha''}}\|_1^2 \Big)^{\frac{1}{2}}&.
\end{align*}
By H\"older's inequality,
\begin{align*}
\int_{\mathfrak{P}(j,t,w,m)}\big|\mathfrak{R}^*f_1 \mathfrak{R}^*f_2 \mathfrak{R}^*f_3\big|\lesssim 2^{\frac{j+r+t}{2}} &\Big(\sum_{\alpha} \| f_{1}\chi_{\tau^j_k \cap t^{j+t}_{w,m}\cap t^{j+2t}_{w^\perp,\alpha}}\|_2^2 \Big)^{\frac{1}{2}}\\		
&\Big(\sum_{\alpha'}\| f_{2}\chi_{\tau^j_{k'} \cap t^{j+t}_{w,m}\cap t^{j+2t}_{w^\perp,\alpha'}}\|_2^2\Big)^{\frac{1}{2}}\\&  \Big(\sum_{w'',\alpha''} \| f_{3}\chi_{\tau^r_{k''} \cap t^{r+t}_{w',m}\cap \varrho^{j+t-r}_{w'',m}\cap t^{2j-r+2t}_{w''^\perp,\alpha''}}\|_2^2 \Big)^{\frac{1}{2}}.
\end{align*}
Thus,
\begin{align*}
\int_{\mathfrak{P}(j,t,w,m)}\big|\mathfrak{R}^*f_1 \mathfrak{R}^*f_2 \mathfrak{R}^*f_3\big|\lesssim 2^{\frac{j+r+t}{2}}\prod_{n=1}^3  \| f_{n}\|_2.
\end{align*}
\begin{flushright}
$\blacksquare$
\end{flushright}

\begin{rema}\label{remarkt}
To see that $\theta^{-\frac{1}{2}}$ is sharp: take $f_1,f_2,f_3$ be $(j,r,t,w,m,\theta)$-triangle type of the form
\begin{align*}
&f_1=\chi_{\tau^j_k \cap t^{j+t}_{w,m}\cap t^{j+2t}_{w^\perp,\alpha}},\\
&f_2=\chi_{\tau^j_{k'} \cap t^{j+t}_{w,m}\cap t^{j+2t}_{w^\perp,\alpha'}}\\
&f_3=\chi_{\tau^r_{k''} \cap t^{r+t}_{w',m}\cap \varrho^{j+t-r}_{w', m}\cap t^{2j-r+2t}_{w'^\perp,\alpha''}}.
\end{align*}
Arguing as in Proposition \ref{Port} and Proposition \ref{Port2}, we see that $\mathfrak{R}^*f_1 \mathfrak{R}^*f_2 \mathfrak{R}^*f_3$ has magnitud $\sim |\tau^j_k \cap t^{j+t}_{w,m}\cap t^{j+2t}_{w^\perp,\alpha}||\tau^j_{k'} \cap t^{j+t}_{w,m}\cap t^{j+2t}_{w^\perp,\alpha'}||\tau^r_{k''} \cap t^{r+t}_{w,m}\cap \varrho^{j+t-r}_{w', m}\cap t^{2j-r+2t}_{w'^\perp,\alpha''}|$ in a large portion of $ \mathfrak{P}(j,t,w,m)$. Thus,
\begin{align*}
\|\mathfrak{R}^*f_1 \mathfrak{R}^*f_2 \mathfrak{R}^*f_3\|&\gtrsim |\mathfrak{P}(j,t,w,m)| 2^{-(j+t+j+2t)} 2^{-(j+t+j+2t)} 2^{-(j+t+2j-r+2t)}\\
&= 2^{-3j} 2^{-4t} 2^{r},
\end{align*}
while we have
\begin{align*}
\| f_{1}\|_2\| f_{2}\|_2\| f_{3}\|_2& \sim 2^{-\frac{1}{2}(j+t+j+2t)} 2^{-\frac{1}{2}(j+t+j+2t)} 2^{-\frac{1}{2}(j+t+2j-r+2t)}\\
&= 2^{-\frac{7}{2}j} 2^{-\frac{9}{2}t} 2^{\frac{1}{2}r}.
\end{align*}

\end{rema}

We introduce some definitions that we use in what follows:
\begin{align*}
&\mathcal{C}_1:=\Big\{c(t^{j+t+\lambda}_{w,\mathfrak{m}_1}\cap t^{j+2t+\lambda}_{w^\perp,\mathfrak{m}_1'}): \quad \tau^j_k \cap t^{j+t}_{w,m}\cap t^{j+t+\lambda}_{w,\mathfrak{m}_1}\cap t^{j+2t+\lambda}_{w^\perp,\mathfrak{m}_1'}\neq \emptyset, \\
&\hspace{46mm}\;\;\;\;\mathfrak{m}_1\in 2^{-(j+t+\lambda)}\mathbb{Z}w^\perp, \mathfrak{m}_1' \in 2^{-(j+2t+\lambda)}\mathbb{Z}w\Big\},\\
&\mathcal{C}_2:=\Big\{c(t^{j+t+\lambda}_{w,\mathfrak{m}_2}\cap t^{j+2t+\lambda}_{w^\perp,\mathfrak{m}_2'}): \quad \tau^j_{k'} \cap t^{j+t}_{w,m}\cap t^{j+t+\lambda}_{w,\mathfrak{m}_2}\cap t^{j+2t+\lambda}_{w^\perp,\mathfrak{m}_2'}\neq \emptyset, \\
&\hspace{46mm}\;\;\;\;\mathfrak{m}_2\in 2^{-(j+t+\lambda)}\mathbb{Z}w^\perp, \mathfrak{m}_2' \in 2^{-(j+2t+\lambda)}\mathbb{Z}w\Big\},\\
&\mathcal{C}_3:=\Big\{c(\varrho^{j+t-r+\lambda}_{\mathfrak{m}_3,m}\cap t^{2j-r+2t+\lambda}_{\mathfrak{m}_3^\perp,\mathfrak{m}_3'}): \quad \tau^r_{k''} \cap t^{r+t}_{w',m}\cap \varrho^{j+t-r+\lambda}_{\mathfrak{m}_3,m}\cap t^{2j-r+2t+\lambda}_{\mathfrak{m}_3^\perp,\mathfrak{m}_3'}\neq \emptyset\\
&\hspace{46mm}\;\;\;\;\mathfrak{m}_3\in \mathbb{S}^1_{j+t-r+\lambda}, \mathfrak{m}_3' \in 2^{-(2j-r+2t+\lambda)}\mathbb{Z}\mathfrak{m}_3\Big\},
\end{align*}

where $c(R)$ is the center of the quadrilateral $R$.

\begin{lema}\label{lemaavO}
For any collection $\{a_{\omega,n}\}_{\omega,n}$ we have
\begin{align*}
\int_{\mathfrak{P}(j,t,w,m)[\lambda]} \big|\prod_{n=1}^3\sum_{\omega\in \mathcal{C}_{n}} a_{\omega,n}e^{2\pi i (\omega,\tfrac{1}{2}|\omega|^2) x}\big|dx\leq C \mathcal{K}(\lambda) \big|{\mathfrak{P}(j,t,w,m)} [\lambda]\big|\prod_{n=1}^3\big(\sum_{\omega\in \mathcal{C}_{n}}  |a_{\omega,n}|^2\big)^{\frac{1}{2}}.
\end{align*}

\end{lema}
\textbf{Proof.}
For the sake of notation compactness, we write
\begin{align*}
&r^{\omega}_1=t^{j+t+\lambda}_{w,\mathfrak{m}_1}\cap t^{j+2t+\lambda}_{w^\perp,\mathfrak{m}_1'},\quad \mbox{with}\quad c(t^{j+t+\lambda}_{w,\mathfrak{m}_1}\cap t^{j+2t+\lambda}_{w^\perp,\mathfrak{m}_1'})=\omega\\
&r^{\omega}_2=t^{j+t+\lambda}_{w,\mathfrak{m}_2}\cap t^{j+2t+\lambda}_{w^\perp,\mathfrak{m}_2'},\quad \mbox{with}\quad c(t^{j+t+\lambda}_{w,\mathfrak{m}_2}\cap t^{j+2t+\lambda}_{w^\perp,\mathfrak{m}_2'})=\omega\\
&r^{\omega}_3=\varrho^{j+t-r+\lambda}_{\mathfrak{m}_3,m}\cap t^{2j-r+2t+\lambda}_{\mathfrak{m}_3^\perp,\mathfrak{m}_3'},\quad \mbox{with}\quad c(\varrho^{j+t-r+\lambda}_{\mathfrak{m}_3,m}\cap t^{2j-r+2t+\lambda}_{\mathfrak{m}_3^\perp,\mathfrak{m}_3'})=\omega.
\end{align*}
We rewrite for $n=1,2,3$,
\begin{align*}
\sum_{\omega\in \mathcal{C}_{n}} a_{\omega,n}e^{2\pi i (\omega,\tfrac{1}{2}|\omega|^2) x}= C\sum_{\omega\in \mathcal{C}_{n}} \int_{r^{\omega}_n}  \big|r^{\omega}_n\big |^{-1}  a_{\omega,n}e^{2\pi i ((\omega,\tfrac{1}{2}|\omega|^2)-(\eta,\tfrac{1}{2}|\eta|^2)) x}e^{2\pi i (\eta,\tfrac{1}{2}|\eta|^2)\cdot x}d\eta.
\end{align*}
We use the Taylor expansion of the exponential
\begin{align*}
e^{2\pi i ((\omega,\tfrac{1}{2}|\omega|^2-(\eta,\tfrac{1}{2}|\eta|^2)) x}=\sum_\gamma c_\gamma ((\omega,\tfrac{1}{2}|\omega|^2)-(\eta,\tfrac{1}{2}|\eta|^2))^\gamma x^\gamma,
\end{align*}
where $\gamma$ denotes a multiindex and the coefficients $c_\gamma$ are decreasing faster than any exponential when $|\gamma|\to\infty$.
Thus,
\begin{align*}
&\sum_{\omega\in \mathcal{C}_{n}} a_{\omega,n}e^{2\pi i (\omega,\tfrac{1}{2}|\omega|^2) x}\\
&= C\sum_\gamma c_\gamma \sum_{\omega\in \mathcal{C}_{n}} \int_{r^{\omega}_n}  \big|r^{\omega}_n\big |^{-1}  a_{\omega,n}\,((\omega,\tfrac{1}{2}|\omega|^2)-(\eta,\tfrac{1}{2}|\eta|^2))^{\gamma} e^{2\pi i(\eta,\tfrac{1}{2}|\eta|^2)\cdot x}d\eta \; x^\gamma.
\end{align*}


Arguing as in Proposition \ref{Port}, we see that the sets $\tilde{r}^{\omega}_n$ are contained in a mild dilation of some $\Big(\big(-(j+2t+\lambda)$, $-(j+t+\lambda)$, $-2(j+t)-\lambda\big)$, $\big((w,0),$ $\big((w,0)\times(m,-1)\big)$, $(m,-1)  \big)\Big)$ parallelepiped, which is precisely the dual parallelepiped of $\mathfrak{P}(j,t,w,m)[\lambda]$. Therefore, we can find an affine transformation $A_{\mathfrak{P}(j,t,w,m)[\lambda]}$, such that 
\begin{align*}
&|(A_{\mathfrak{P}(j,t,w,m)[\lambda]}((\omega,\tfrac{1}{2}|\omega|^2)-(\eta,\tfrac{1}{2}|\eta|^2)))^{\gamma}|\lesssim 1,\;\;\eta\in r^{\omega}_n\\
&|((A_{\mathfrak{P}(j,t,w,m)[\lambda]})^{-1}(x))^\gamma|\lesssim 1,\;\;  x\in {\mathfrak{P}(j,t,w,m)[\lambda]}.
\end{align*}
Thus,
\begin{align*}
&\Big|\sum_{\omega\in \mathcal{C}_{n}} a_{\omega,n}e^{2\pi i (\omega,\tfrac{1}{2}|\omega|^2) x}\Big|= \Big|\sum_\gamma c_\gamma \mathfrak{R}^*\phi_{\gamma,t}(x)((A_{\mathfrak{P}(j,t,w,m)[\lambda]})^{-1}(x))^\gamma
\Big|\leq C \sum_\gamma |c_\gamma|  | \mathfrak{R}^*\phi_{\gamma,n}(x)|
\end{align*}
where
\begin{align*}
\phi_{\gamma,n}(\eta)= \sum_{\omega\in \mathcal{C}_{n}} \chi_{r^{\omega}_n}  \big|r^{\omega}_n\big |^{-1} (A_{\mathfrak{P}(j,t,w,m)[\lambda]}((\omega,\tfrac{1}{2}|\omega|^2)-(\eta,\tfrac{1}{2}|\eta|^2)))^{\gamma} a_{\omega,n}.
\end{align*}
Integrating over $\mathfrak{P}(j,t,w,m)[\lambda]$,
\begin{align*}
\int_{\mathfrak{P}(j,t,w,m)[\lambda]} \big|\prod_{n=1}^3\sum_{\omega\in \mathcal{C}_{n}} a_{\omega,n}e^{2\pi i (\omega,\tfrac{1}{2}|\omega|^2) x}\big|dx\lesssim  \sum_{\gamma_1,\gamma_2,\gamma_3} |c_{\gamma_1}||c_{\gamma_2}||c_{\gamma_3}|  \int_{\mathfrak{P}(j,t,w,m)[\lambda]} \prod_{n=1}^3|\mathfrak{R}^*\phi_{\gamma,n}|dx.
\end{align*}
Thus,
\begin{align*}
&\int_{\mathfrak{P}(j,t,w,m)[\lambda]} \big|\prod_{n=1}^3\sum_{\omega\in \mathcal{C}_{n}} a_{\omega,n}e^{2\pi i (\omega,\tfrac{1}{2}|\omega|^2) x}\big|dx\\
&\lesssim  \sum_{\gamma_1,\gamma_2,\gamma_3} |c_{\gamma_1}||c_{\gamma_2}||c_{\gamma_3}| 2^{\frac{j+r+t}{2}} \mathcal{K}(\lambda) \prod_{n=1}^3 \|\phi_{\gamma,n}\|_{L^2}\\
&=\sum_{\gamma_1,\gamma_2,\gamma_3} |c_{\gamma_1}||c_{\gamma_2}||c_{\gamma_3}| 2^{\frac{j+r+t}{2}} \mathcal{K}(\lambda)\\
&\hspace{35mm}\prod_{n=1}^3 \Big(\int \big| \sum_{\omega\in \mathcal{C}_{n}} \chi_{r^{\omega}_n}  \big|r^{\omega}_n\big |^{-1} (A_{\mathfrak{P}(j,t,w,m)[\lambda]}((\omega,\tfrac{1}{2}|\omega|^2)-(\eta,\tfrac{1}{2}|\eta|^2)))^{\gamma} a_{\omega,n}\big|^2\Big)^{\frac{1}{2}}\\
&\lesssim\sum_{\gamma_1,\gamma_2,\gamma_3} |c_{\gamma_1}||c_{\gamma_2}||c_{\gamma_3}| 2^{\frac{j+r+t}{2}}\mathcal{K}(\lambda)\prod_{n=1}^3   |r^{\omega}_n|^{-\frac{1}{2}}\Big(\sum_{\omega\in \mathcal{C}_{n}} |a_{\omega,n}|^2 \Big)^{\frac{1}{2}}\\
&\sim\sum_{\gamma_1,\gamma_2,\gamma_3} |c_{\gamma_1}||c_{\gamma_2}||c_{\gamma_3}| |\mathfrak{P}(j,t,w,m)[\lambda]|\mathcal{K}(\lambda)\prod_{n=1}^3 \Big(\sum_{\omega\in \mathcal{C}_{n}} |a_{\omega,n}|^2 \Big)^{\frac{1}{2}},
\end{align*}
and by the decay of the $c_{\gamma_t}$ we conclude the result.
\begin{flushright}
$\blacksquare$
\end{flushright}

We will use the well known estimate
\begin{proposition}\label{unmedio}
$$\|\mathfrak{R}^* f\|_{L^2(B_R)}\lesssim R^{\frac{1}{2}}\|f\|_2,\quad \mbox{for every}\; R\geq 1. $$
\end{proposition}

\textbf{Proof of Proposition \ref{P1}.}

Fix $z\in \mathfrak{P}(j,t,w,m)[2\lambda]$. For $x,y\in \mathfrak{P}(j,t,w,m)[\lambda](z)$ we have
\begin{align*}
\sum_{\omega\in \mathcal{C}_n}&\mathfrak{R}^*f_n \chi_{r^{\omega}_n}(x)\\
&= \sum_{\omega\in \mathcal{C}_n} \int_{r^{\omega}_n} f_n (\eta) e^{2\pi i ((\omega,\tfrac{1}{2}|\omega|^2)-(\eta,\tfrac{1}{2}|\eta|^2))\cdot(x-y)} e^{-2\pi i (\eta,\tfrac{1}{2}|\eta|^2) y}d\sigma(\eta) e^{-2\pi i (\omega,\tfrac{1}{2}|\omega|^2) (x-y)}.
\end{align*}
Defining $f_{n,\gamma,\omega}(\eta)=f_n(\eta) \chi_{r^{\omega}_n} (A_{\mathfrak{P}(j,t,w,m)[\lambda]}((\omega,\tfrac{1}{2}|\omega|^2)-(\eta,\tfrac{1}{2}|\eta|^2)))^{\gamma} $, we have, as in the proof of the previous lemma,
\begin{align*}
\sum_{\omega\in \mathcal{C}_n}\mathfrak{R}^*f_n \chi_{r^{\omega}_n}(x)= \sum_{\gamma_n} c_{\gamma_n}  ((A_{\mathfrak{P}(j,t,w,m)[\lambda]})^{-1}(x-y))^\gamma\sum_{\omega\in \mathcal{C}_n}\mathfrak{R}^*f_{n,\gamma,\omega}(y)e^{-2\pi i (\omega,\tfrac{1}{2}|\omega|^2) (x-y)}.
\end{align*}
Thus,
\begin{align*}
&\Big\|\prod_{n=1}^3\sum_{\omega\in \mathcal{C}_n} \mathfrak{R}^*f_n \chi_{r^{\omega}_n}\Big\|_{L^1 (\mathfrak{P}(j,t,w,m)[\lambda](z))}\\
&\hspace{15mm}\lesssim \sum_{\gamma_1} c_{\gamma_1}  \sum_{\gamma_2} c_{\gamma_2} \sum_{\gamma_3} c_{\gamma_3} \Big\|\prod_{n=1}^3\sum_{\omega\in \mathcal{C}_n}\mathfrak{R}^*f_{n,\gamma,\omega}(y)e^{-2\pi i (\omega,\tfrac{1}{2}|\omega|^2) (\cdot-y)}\Big\|_{L^1 (\mathfrak{P}(j,t,w,m)[\lambda](z))}.
\end{align*}
Using Lemma \ref{lemaavO} and averaging with respect to $y\in \mathfrak{P}(j,t,w,m)[\lambda](z)$
\begin{align*}
&\Big\|\prod_{n=1}^3\sum_{\omega\in \mathcal{C}_n}\mathfrak{R}^*f_n \chi_{r^{\omega}_n}\Big\|_{L^1 (\mathfrak{P}(j,t,w,m)[\lambda](z))}\\
&\hspace{15mm}\lesssim \sum_{\gamma_1} c_{\gamma_1}  \sum_{\gamma_2} c_{\gamma_2} \sum_{\gamma_3} c_{\gamma_3} \mathcal{K}(\lambda)\Big\|\prod_{n=1}^3\Big(\sum_{\omega\in \mathcal{C}_n}\big|\mathfrak{R}^*f_{n,\gamma,\omega}\big|^2\Big)^{\frac{1}{2}}\Big\|_{L^1 (\mathfrak{P}(j,t,w,m)[\lambda](z))}.
\end{align*}
As we can find a collection $$\mathfrak{P}(j,t,w,m)[2\lambda]\subset \bigcup_z \mathfrak{P}(j,t,w,m)[\lambda](z),$$
summing on $z$, we get
\begin{align*}
&\Big\|\prod_{n=1}^3\sum_{\omega\in \mathcal{C}_n}\mathfrak{R}^*f_n \chi_{r^{\omega}_n}\Big\|_{L^1 (\mathfrak{P}(j,t,w,m)[2\lambda])}\\
&\hspace{15mm}\lesssim \sum_{\gamma_1} c_{\gamma_1}  \sum_{\gamma_2} c_{\gamma_2} \sum_{\gamma_3} c_{\gamma_3}  \mathcal{K}(\lambda)\Big\|\prod_{n=1}^3\Big(\sum_{\omega\in \mathcal{C}_n}\big|\mathfrak{R}^*f_{n,\gamma,\omega}\big|^2\Big)^{\frac{1}{2}}\Big\|_{L^1 (\mathfrak{P}(j,t,w,m)[2\lambda])}.
\end{align*}
Now, as $\mathfrak{P}(j,t,w,m)[2\lambda]\subset B_{2\lambda+2j+2t}$, we get \begin{align*}
&\chi_{\mathfrak{P}(j,t,w,m)[2\lambda]}\Big(\sum_{\omega\in \mathcal{C}_n}\big|\mathfrak{R}^*f_{n,\gamma,\omega}\big|^2\Big)^{\frac{1}{2}}\lesssim \phi_{B_{2\lambda+2j+2t}}\Big(\sum_{\omega\in \mathcal{C}_n}\big|\mathfrak{R}^* f_{n,\gamma,\omega}\chi_{r^\omega_{n}}\big|^2 \Big)^{\frac{1}{2}}\\& \lesssim \Big(\sum_{\overline{k}\in 2^{-(j+t+\lambda)}\mathbb{Z}^2} \sum_{\omega\in \mathcal{C}_n:  \;r^\omega_{n}\cap \tau^{j+t+\lambda}_{\overline{k}}\neq \emptyset}\big| \phi_{B_{2\lambda+2j+2t}}\mathfrak{R}^* f_{n,\gamma,\omega}\chi_{r^\omega_{n}\cap\tau^{j+t+\lambda}_{\overline{k}}}\big|^2\Big)^{\frac{1}{2}}.
 \end{align*}

Let $\psi_{\tau^{j+t+\lambda}_{\overline{k}}}$ be a Schwartz function which is comparable to $1$ on the parallelepiped containing $\tilde{\tau}^{j+t+\lambda}_{\overline{k}}$, and whose compactly supported Fourier transform satisfies
$$|\widehat{\psi_{\tau^{j+t+\lambda}_{\overline{k}}}}|\lesssim |\accentset{\circ}{\tau}^{j+t+\lambda}_{\overline{k}}|^{-1} \chi_{\accentset{\circ}{\tau}^{j+t+\lambda}_{\overline{k}}}.$$ 
We have by Jensen inequality that
$$\big|\phi_{B_{2\lambda+2j+2t}} \mathfrak{R}^* f_{n,\gamma,\omega}\chi_{r^\omega_{n}\cap\tau^{j+t+\lambda}_{\overline{k}}}\big|^2\lesssim \big|\phi_{B_{2\lambda+2j+2t}} \mathfrak{R}^* f_{n,\gamma,\omega}\chi_{r^\omega_{n}\cap\tau^{j+t+\lambda}_{\overline{k}}}*\varphi_{\overline{k}}\big|^2*|\widehat{\psi_{\tau^{j+t+\lambda}_{\overline{k}}}}|,$$ where $\varphi_{\overline{k}}$ satisfies $\widehat{\varphi}_{\overline{k}}\Big(\phi_{B_{2\lambda+2j+2t}} \mathfrak{R}^* f_{n,\gamma,\omega}\chi_{r^\omega_{n}\cap\tau^{j+t+\lambda}_{\overline{k}}}\Big)^\vee\sim \Big(\phi_{B_{2\lambda+2j+2t}} \mathfrak{R}^* f_{n,\gamma,\omega}\chi_{r^\omega_{n}\cap\tau^{j+t+\lambda}_{\overline{k}}}\Big)^\vee$. 

Hence, 
\begin{align*}
&\chi_{\mathfrak{P}(j,t,w,m)[2\lambda]}\Big(\sum_{\omega\in \mathcal{C}_n}\big|\mathfrak{R}^*f_{n,\gamma,\omega}\big|^2\Big)^{\frac{1}{2}}\lesssim  |\accentset{\circ}{\tau}^{j+t+\lambda}_{\overline{k}}|^{-\frac{1}{2}}\\
&\hspace{5mm}\times\Big(\sum_{\overline{k}\in 2^{-(j+t+\lambda)}\mathbb{Z}^2} \big(\sum_{\omega\in \mathcal{C}_n:  \;r^\omega_{n}\cap \tau^{j+t+\lambda}_{\overline{k}}\neq \emptyset}\big| \phi_{B_{2\lambda+2j+2t}} \mathfrak{R}^* f_{n,\gamma,\omega}\chi_{r^\omega_{n}\cap\tau^{j+t+\lambda}_{\overline{k}}}*\varphi_{\overline{k}}\big|^2\big)*\chi_{\accentset{\circ}{\tau}^{j+t+\lambda}_{\overline{k}}}\Big)^{\frac{1}{2}}.
 \end{align*}
 Thus,
\begin{align*}
&\Big\|\prod_{n=1}^3\sum_{\omega\in \mathcal{C}_n}\mathfrak{R}^*f_n \chi_{r^{\omega}_n}\Big\|_{L^1 (\mathfrak{P}(j,t,w,m)[2\lambda])}\\
&\lesssim \sum_{\gamma_1} c_{\gamma_1}  \sum_{\gamma_2} c_{\gamma_2} \sum_{\gamma_3} c_{\gamma_3}\mathcal{K}(\lambda) |\accentset{\circ}{\tau}^{j+t+\lambda}|^{-\frac{3}{2}}  \\
&\hspace{5mm}\Big\|\prod_{n=1}^3\Big(\sum_{\overline{k}\in 2^{-(j+t+\lambda)}\mathbb{Z}^2} \big( \sum_{\omega\in \mathcal{C}_n:  \;r^\omega_{n}\cap \tau^{j+t+\lambda}_{\overline{k}}\neq\emptyset}\big| \phi_{B_{2\lambda+2j+2t}} \mathfrak{R}^* f_{n,\gamma,\omega}\chi_{r^\omega_{n}\cap\tau^{j+t+\lambda}_{\overline{k}}}*\varphi_{\overline{k}}\big|^2\big)*\chi_{\accentset{\circ}{\tau}^{j+t+\lambda}_{\overline{k}}}\Big)^{\frac{1}{2}}\Big\|_{L^1}.
\end{align*}
Using Theorem \ref{GuthT} and Placherel's theorem, 
\begin{align*}
&\Big\|\prod_{n=1}^3\sum_{\omega\in \mathcal{C}_n}\mathfrak{R}^*f_n \chi_{r^{\omega}_n}\Big\|_{L^1 (\mathfrak{P}(j,t,w,m)[2\lambda])}\\
&\lesssim \sum_{\gamma_1} c_{\gamma_1}  \sum_{\gamma_2} c_{\gamma_2} \sum_{\gamma_3} c_{\gamma_3} \mathcal{K}(\lambda)|\accentset{\circ}{\tau}^{j+t+\lambda}|^{-\frac{3}{2}} 2^{\frac{j+r+t}{2}} 2^{3(j+t+\lambda)}\\
&\hspace{20mm}\prod_{n=1}^3\Big( \sum_{\overline{k}\in 2^{-(j+t+\lambda)}\mathbb{Z}^2}\Big\|\sum_{\omega\in \mathcal{C}_n:  \;r^\omega_{n}\cap \tau^{j+t+\lambda}_{\overline{k}}\neq\emptyset}\big| \phi_{B_{2\lambda+2j+2t}} \mathfrak{R}^* f_{n,\gamma,\omega}\chi_{r^\omega_{n}\cap\tau^{j+t+\lambda}_{\overline{k}}}*\varphi_{\overline{k}}\big|^2\Big\|_{L^1}\Big)^{\frac{1}{2}}\\
&\sim \sum_{\gamma_1} c_{\gamma_1}  \sum_{\gamma_2} c_{\gamma_2} \sum_{\gamma_3} c_{\gamma_3}\mathcal{K}(\lambda) |\accentset{\circ}{\tau}^{j+t+\lambda}|^{-\frac{3}{2}} 
2^{\frac{j+r+t}{2}} 2^{3(j+t+\lambda)}\prod_{n=1}^3\Big( \sum_{\omega\in \mathcal{C}_n}\Big\|  \phi_{B_{2\lambda+2j+2t}} \mathfrak{R}^* f_{n,\gamma,\omega}\Big\|^2_{L^2}\Big)^{\frac{1}{2}}.
\end{align*}
Using Proposition \ref{unmedio}, and the rapidly decreasing of $ \phi_{B_{2\lambda+2j+2t}}$ away from $B_{2\lambda+2j+2t}$,
\begin{align*}
&\Big\|\prod_{n=1}^3\sum_{\omega\in \mathcal{C}_n}\mathfrak{R}^*f_n \chi_{r^{\omega}_n}\Big\|_{L^1 (\mathfrak{P}(j,t,w,m)[2\lambda])}\\&\lesssim  \sum_{\gamma_1} c_{\gamma_1}  \sum_{\gamma_2} c_{\gamma_2} \sum_{\gamma_3} c_{\gamma_3}\mathcal{K}(\lambda) |\accentset{\circ}{\tau}^{j+t+\lambda}|^{-\frac{3}{2}} 
2^{\frac{j+r+t}{2}} 2^{6(j+t+\lambda)}\prod_{n=1}^3\Big( \sum_{\omega\in \mathcal{C}_n}\Big\|  f_{n,\gamma,\omega}\Big\|^2_{L^2}\Big)^{\frac{1}{2}}.
\end{align*}
Using that $f_{n,\gamma,\omega}\lesssim f_n\chi_{r^\omega_{n}}$, and the decay of $c_{\gamma_1}, c_{\gamma_2}, c_{\gamma_3}$,
\begin{align*}
\Big\|\prod_{n=1}^3\sum_{\omega\in \mathcal{C}_n}\mathfrak{R}^*f_n \chi_{r^{\omega}_n}\Big\|_{L^1 (\mathfrak{P}(j,t,w,m)[2\lambda])}&\lesssim \mathcal{K}(\lambda) |\accentset{\circ}{\tau}^{j+t+\lambda}|^{-\frac{3}{2}} 
2^{\frac{j+r+t}{2}} 2^{6(j+t+\lambda)}\prod_{n=1}^3\Big( \sum_{\omega\in \mathcal{C}_n}\Big\|  f_n\chi_{r^\omega_{n}}\Big\|^2_{L^2}\Big)^{\frac{1}{2}}\\
&= \mathcal{K}(\lambda) |\accentset{\circ}{\tau}^{j+t+\lambda}|^{-\frac{3}{2}} 
2^{\frac{j+r+t}{2}} 2^{6(j+t+\lambda)}\prod_{n=1}^3\|  f_{n}\|_{L^2}\\
&\sim   \mathcal{K}(\lambda) 2^{\frac{j+r+t}{2}} \prod_{n=1}^3\|  f_{n}\|_{L^2}.
\end{align*}

\begin{flushright}
$\blacksquare$
\end{flushright}

\textbf{Acknowledgments:} I want to thank Keith Rogers and Ana Vargas for their support. I learned about the restriction problem from them during my PhD years, and they encouraged me to work on it. Thanks to Sanghyuk Lee for his helpful conversation regarding the main theorem. Thanks to Felipe Ponce for a careful reading and for pointing out some mistakes. Thanks to Felipe Linares and Emanuel Carneiro for these past years at IMPA. 


\end{document}